\documentclass[11pt]{article} 
\usepackage{amsmath,amssymb}
\usepackage[utf8]{inputenc}

\usepackage{geometry} 
\usepackage{graphicx} 

\usepackage{booktabs} 
\usepackage{array} 
\usepackage{paralist} 
\usepackage{verbatim} 
\usepackage{subfig} 

\usepackage{fancyhdr} 
\usepackage{sectsty}
\allsectionsfont{\sffamily\mdseries\upshape} 

\usepackage[nottoc,notlof,notlot]{tocbibind} 
\usepackage[titles,subfigure]{tocloft} 

\title{Linear algebra and quantum algorithm\footnote{This paper wirtten by the finantial support of Brain Korea 21}}
\author{BongJu Kim}
\date{} 

\begin{document}
\maketitle
\abstract{We introduce quantum algorithm and the mathematical structure of quantum computer. Quantum algorithm is expressed by linear algebra on a finite dimensional complex inner product space. The mathematical formulations of quantum mechanics had been established in around 1930, by von Neumann. The formulation uses functional analysis, linear algebra and probability theory. The knowledge of mathematical formulations of QM is enough quantum mechanical knowledge for approaching to quantum algorithm and it might be efficient way for mathematicians that starting with the mathematical formulations of QM. We explain the mathematical formulations of quantum mechanics briefly, quantum bits, quantum gates, quantum discrete Fourier transformation, Deutsch's algorithm and Shor's algorithm.}
\section{Introduction}

As quantum computer hardware production, which seemed a long way off, has made some progresses recently, much attention is also being paid to the study of quantum algorithm. The class of decision problems which solvable by a quantum computer in polynomial time is called BQP(bounded error quantum polynomial time). Although BQP is not perfectly identified yet, It was proved that many important and hard decision problems belong to BQP. Cryptologists also regard quantum computing as a realizable threat. For example, Shor's algorithm can broke a cypher which relying on the difficulty of the discrete logarithm such as RSA or ECC(elliptic curve cryptography). Cryptologists are preparing for the quantum computing era. This research field called post-quantum cryptography.

Quantum algorithm is expressed by linear algebra on a finite dimensional complex inner product space.   The part that need to know about QM is just the mathematical formulation of quantum mechanics which is formulated by probability theory, linear algebra and functional analysis.  So, quantum algorithm is just a mathematical problem. In fact, many mathematicians, such as Peter Shor\footnote{His prime factorization quantum algorithm made a sensational impact and triggered many researches about quantum computing and financial  investments because it can broke a strong cryptography system.}, Michael Freedman\footnote{A fields medal winner mathematician. He works in  Microsoft Quantum – Santa Barbara.}, research quantum algorithm.

This paper is the lecture note that I wrote for the (about) six-hour lecture that I spoke in Quantum algorithm seminar during 2019 spring semester. I don't know about and unfamiliar with physics\footnote{I do not interested in sciences but because of Shor's algorithm, I had become interested in the mathematical formulations of quantum mechanics. Fortunately, }, but it did not take long time to approach quantum algorithm. I expect that the readers will be able to understand it easily.

\section{The mathematical formulations of quantum mechanics}

Before quantum mechanics, one of the main purpose in physics was to find `` the trajectory of a particle'', $x:(a,b)\rightarrow \mathbb R^3 $ mathematically, from initial location and momentum of the particle and mechanical principles which are mathematically formulated mainly in a system of partial differential equations. This way had been established after 17 century-the birth of physics. It was believed that the initial locations and momentums of a physical system determines perfectly the future of the physical system. There was also an extreme argument, in this way, known as ``Laplace's demon'' by a French mathematician Pierre Simon Laplace.
 In this direction, Newtonian mechanics, Lagrangian mechanics, Hamiltonian mechanics\footnote{However, Lagrangian mechanics and Hamiltonian mechanics seems like prepare quantum mechanics and quantum field theory.} and the theory of relativity were very successful in the description of the macroscopic physical world. 

However, in atomic scale  (about $10^{-9}$m) physics, finding `` the trajectory of a particle''  is an unattainable purpose according to quantum physics. In atomic scale physics, one can't know what physical event will be happened but only can say about the ``distribution of probability''. Also, `` the trajectory of a particle'' does not make sense in this scale. 

Suppose that you want to know about the momentum of a particle with mass $m$ in a specific potential environment\footnote{Let's consider 1-dimensional case.} described as the real-valued function $V(x,t)$. Then you solve the Schrodinger's equation \[i\hbar\frac{\partial}{\partial t}\psi(x,t) = -\frac{\hbar^2}{2m}\frac{\partial^2}{\partial x^2}\psi(x,t) +V(x,t)\psi(x,t)\] \footnote{where $\hbar$ is the Dirac's constant $1.054571817\times 10^{-34}$J$\cdot$s} to find ``the wave function of the particle'' $\psi(x,t)$ which have all quantum mechanical information about the physical system.  A  wave function is a complex valued function and an element of a complex Hilbert space $\mathcal H$(complete inner product space over $\mathbb C$, $L^2$ space usually). As you know, any constant multiple of a solution of the linear partial differential equation is also a solution\footnote{However, the zero function do not fit to describe a physical system. So, we consider only non-zero complex functions.}, we take the solution with unit norm. Now, operate ``the momentum operator'' \[\hat{p}:=  \frac{\hbar }{i}\frac{ \partial}{\partial x}\] which is an Hermitian on the wave equation $\psi(x,t)$ and compute the inner product of $\psi(x,t)$ and $\hat{p}\psi(x,t)$:\[\int\psi(x,t)(\hat{p}\psi(x,t))dx=\int\psi(x,t)\frac{\hbar }{i}\frac{ \partial\psi(x,t)}{\partial x}dx.\] Then this value means that the expectation value of the momentum of the particle. Because any constant multiple of a solution of the linear partial differential equation is also a solution, we give a equivalence relation on $\mathcal H\setminus\{0\}$:\[\text{For any \,}\phi , \psi\in \mathcal H\setminus\{0\},\;\;\phi \sim \psi \text{\;\;iff\;\;}\exists c\in \mathbb C \text{\;\;s.t.\;\;} \phi=c\psi.\]
\\

Generally, quantum mechanics can be mathematically formulated as follows:
\\

1. A quantum mechanical system associated with a separable\footnote{An inner product space is trivially a normed linear space(Banach space). If a Banach space is not separable, then there is no Schauder basis. If $S$ is a Schauder basis of a Banach space $\mathcal X$, Span$S$ is dense subset of $\mathcal X$.} complex Hilbert space $\mathcal H$. A quantum sate is described by a 1-dimensional subspace of $\mathcal H$. Especially. the zero element in $\mathcal H$ is do not fit and a quantum sate is exactly associated with an element  of  complex projective Hilbert space $(\mathcal H\setminus\{0\})/\sim$. Therefore, we can take a element with unit norm as a representative.

2. Let $\mathcal H_1, \mathcal H_2$ describe two quantum mechanical systems respectably. Then, the Hilbert space which describes the composition of two quantum mechanical systems is $\mathcal H_1\otimes \mathcal H_2$.

3. Physical observables \footnote{For example,  position, momentum, energy, spin, etc. } are described by Hermitian operators on $\mathcal H$.

4. The expectation value of an observable $\hat A$ of a quantum mechanical system in the state represented by the unit element $\psi \in \mathcal H$ is the inner product of $\psi$ and $\hat A \psi$.

5. Physical symmetries in qunatum mechanics are represented by unitary or anti-unitary operators.\footnote{Due to Wigner's theorem.}

6. Let an observable represented by $\hat A$ in a quantum mechanical system has a discrete spectrum$\{\lambda_i\;|\;i=1,2,\dots \}$. Then, the result of the experimental measurement is one of the eigenvalues $\lambda_i$ and the probability that we get the result $\lambda_i$ is the inner product of $\psi$ and $\hat P_i\psi$ where $\hat P_i$ is the projection operator corresponding to $\lambda_i$.

\section{Quantum bits}

A quantum bit(qubit) is the unit of information in quantum computing, and one of the unit elements of 2-dimensional complex Hilbert space $H$ with the inner product\[(\cdot,\cdot):\mathcal H\times\mathcal  H\longrightarrow \mathbb C.\]
As a bit can be physically implemented by two different voltage or power on-off,  Quantum bit can be physically implemented by any two-state quantum mechanical system such as two states of spin of an electron or two states of polarization of a photon.

1-qubit with an orthonormal basis $\{b_0, b_1\}$ represented by \[u=c_0b_0+c_1b_1 \in \mathcal H \;\;\;\text{where}\;\;(u,u)=c_0\bar c_0+c_1\bar c_1=|c_0|^2+|c_1|^2=1\]and \[u=c_0b_0+c_1b_1=c_0\begin{bmatrix}1 \\ 0 \end{bmatrix}+c_1\begin{bmatrix} 0 \\ 1 \end{bmatrix}=\begin{bmatrix} c_0 \\ c_1 \end{bmatrix}.\] $b_0,b_1$ means the bit 0 and 1. The measurement of qubit is probabilistic. The sample space of 1-qubit measurement is  $\{b_0, b_1\}$ and the probability of event $b_i$ is $|c_i|^2$. 

$n$-qubit system associated with $\mathcal H^{\otimes n}$ so that represented by a unit element of a $2^n$-dimensional complex Hilbert space\footnote{Trivially, any higher $m$-dimensional complex Hilbert space and $m^n$-dimensional complex Hilbert space are also possible. } with an orthonormal basis $\{b_i\;|\;i=0,1,\dots ,2^n-1\}$: \[v=\sum_{i=1}^{2^n}c_ib_i\in \mathcal H^{\otimes n}\;\;\;\text{where}\;\;\sum_{i=1}^{2^n}|c_i|^2=1.\]

Here, the $2^n$-dimensional complex projective Hilbert space is the stage for quantum algorithms are performed. It is easy to see that $b_i\rightarrow i$ denotes all possible bit from $0$ to $2^n-1$ i.e. the basis $\{b_i\;|\;i=0,1,\dots ,2^n-1\}$ is the sample space of the n-qubit measurement and $P(b_i)=c_i\bar c_i=|c_i|^2$.

Trivially, there is an element $u$ in $\mathcal H\otimes \mathcal H$ such that $u$ is not a Kronecker product  $a\otimes b$ where $a,b\in \mathcal H$. It is the mathematical formulation of quantum entanglement. For example, 2-qubit system is generally \[c_0(b_0\otimes b_0)+c_1(b_0\otimes b_1)+c_2(b_1\otimes b_0)+c_3(b_1\otimes b_1)=\begin{bmatrix} c_0 \\ c_1\\ c_2 \\ c_3\end{bmatrix}\] and the Kronecker product of $u=u_0b_0+u_1b_1$ and $v=v_0b_0+v_1b_1$ is \[u\otimes v= u_0v_0(b_0\otimes b_0)+u_1v_0(b_0\otimes b_1)+u_0v_1(b_1\otimes b_0)+u_1v_1(b_1\otimes b_1)=\begin{bmatrix} u_0v_0 \\ u_1v_0\\ u_0v_1 \\ u_1v_1\end{bmatrix}.\] Therefore, there are many element which is not a Kronecker product of two elements in $\mathcal H$ such as \[w=\frac{1}{\sqrt 2}(b_1\otimes b_0)+\frac{1}{\sqrt 2}(b_0\otimes b_1)=\frac{1}{\sqrt 2}\begin{bmatrix} 0 \\ 1\\ 1 \\ 0\end{bmatrix}.\]$w$ is not a a Kronecker product of two elements in $\mathcal H$ since $\mathbb C$ has no zero divisor. In quantum algorithm, two qubits can be entangled as a result of a quantum gate operation.

\section{Quantum gates and what is a quantum algorithm}
A quantum gate on $n$-qubit is a unitary linear map $U$ on $2^n$-dimensional complex Hilbert space $ \mathcal H^{\otimes n}$ and represented by $2^n\times 2^n$ unitary matrix \footnote{Recall that if a matrix $U$ satisfies $UU^*=U^*U=I$, $U$ is unitary }. Since a unitary map preserves the norm of elements, the result $Uv$ is also unit. A quantum gate changes the probability distribution on the basis. Suppose that there is a  problem and we prepared enough qubits to express the answer of the problem. This means that the set basis $B$ of $ \mathcal H^{\otimes n}$ contains the answer. Now, \textit{a quantum algorithm to solve the problem is a sequence of quantum gates which makes the probability of the answer of the problem, denoted by a basis element $b^* \in  B$, higher enough so that we can get the answer quickly by iterating performance of the quantum algorithm.} For example, suppose that a quantum algorithm have the probability of the answer is $1/5$ then, the probability that the results of 15 performances never meet the answer is $(4/5)^{15}\approx 0.03$. Trivially, since the result of quantum algorithm is probabilistic, we should verify whether the result is really the answer or not. We can use a classical computer to check it.

Since quantum gate is unitary, it's invertible. Therefore, a quantum computation can be traced back from the result, and it preserves all informations. this is one point that quantum computations differ from classical computations.

Following matrices(quantum gates) act on a single qubit. The Hadamard matrix(gate) is
\[H := \frac{1}{\sqrt{2}} \begin{bmatrix} 1 & 1 \\ 1 & -1 \end{bmatrix}.\] Observe that 
\[Hu=\frac{1}{\sqrt{2}} \begin{bmatrix} 1 & 1 \\ 1 & -1 \end{bmatrix} \begin{bmatrix} c_0 \\ c_1 \end{bmatrix}=\frac{1}{\sqrt{2}}\begin{bmatrix} c_0+c_1 \\ c_0-c_1 \end{bmatrix}\]and it makes a superposition if $u$ is a basis bit $b_0$ or $b_1$. The $X$-gate is
\[X := \begin{bmatrix} 0 & 1 \\ 1 & 0 \end{bmatrix}.\] It changes the coefficients of a qubit:\[Xu=\begin{bmatrix} 0 & 1 \\ 1 & 0 \end{bmatrix} \begin{bmatrix} c_0 \\ c_1 \end{bmatrix}=\begin{bmatrix} c_1 \\ c_0 \end{bmatrix}.\] It is analogous to the classical NOT gate since it flips the bit when it acts on a basis bit. The twist gates

\[T(\alpha) := \begin{bmatrix} 1 & 0 \\ 0 & e^{i\alpha} \end{bmatrix}\] do not change the probability distribution but   change the argument of a coefficient:
\[T(\alpha)u = \begin{bmatrix} 1 & 0 \\ 0 & e^{i\alpha} \end{bmatrix} \begin{bmatrix} c_0 \\ c_1 \end{bmatrix}=\begin{bmatrix} c_0 \\ e^{i\alpha}c_1 \end{bmatrix}.\]

There are many matrices(quantum gates) act on two qubit. But here, we present a very important quantum gate which involves a quantum entanglement. The quantum gate is CNOT(controlled-not) gate
\[\wedge_1(X):=\begin{bmatrix} 
1 & 0 & 0 & 0 \\
0 & 1 & 0 & 0 \\
0 & 0 & 0 & 1\\
0 & 0 & 1 & 0\\

\end{bmatrix}.\] It acts as identity gate for the first qubit and as X-gate(which is analogous to the classical NOT gate). Observe that  \[\wedge_1(X)v:=\begin{bmatrix} 
1 & 0 & 0 & 0 \\
0 & 1 & 0 & 0 \\
0 & 0 & 0 & 1\\
0 & 0 & 1 & 0\\

\end{bmatrix}\begin{bmatrix} v_0 \\ v_1 \\ v_2 \\ v_3 \end{bmatrix}=\begin{bmatrix} v_0 \\ v_1 \\ v_3 \\ v_2 \end{bmatrix}.\] Especially, 
\begin{align*}
\wedge_1(X)(b_0\otimes b_0)&=b_0\otimes b_0,\\
\wedge_1(X)(b_1\otimes b_0)&=b_1\otimes b_0,\\
\wedge_1(X)(b_0\otimes b_1)&=b_1\otimes b_1,\\
\wedge_1(X)(b_1\otimes b_1)&=b_0\otimes b_1.\\
\end{align*}
i.e.
\[\wedge_1(X)(b_j\otimes b_i)=b_{j\oplus i}\otimes b_i.\] 
where $\oplus$ is the addition in $\mathbb Z_2$. It is showed that  CNOT gate is enough for any quantum circuit involving a quantum entanglement and we do not need any other quantum entanglement-involving gates.

\section{Quantum discrete Fourier transformation}
Quantum discrete Fourier transformation is an important transformation in many quantum algorithms. This is just  discrete Fourier transformation \[\hat f(k)=\frac{1}{\sqrt N}\sum_{j=0}^{N-1}e^{2\pi ijk/N}f(j)\] on qubits. For a basis element $b_k$ of $\mathcal H^{\otimes n}$, the quantum discrete Fourier transformation $\mathcal F_n$ on $n$-qubit is \[\mathcal F_n(b_k)=\frac{1}{\sqrt {2^n}}\sum_{j=0}^{2^n-1}e^{2\pi ijk/2^n}b_j.\] As you know, if $N | k$, then \[\frac{1}{\sqrt N}\sum_{j=0}^{N-1}e^{2\pi ijk/N}=1\] and it is 0 if $N$ not divide $k$. Quantum discrete Fourier transformation is due to a mathematician and a cryptographer Don Coppersmith.\footnote{Coppersmith, D. An approximate Fourier transform useful in quantum factoring. Technical Report RC19642, IBM. 1994}

Let us denote $e^{2\pi i/2^n}$ by $\zeta_{2^n}$. Quantum discrete Fourier transformation on $n$-qubit represented by the unitary matrix \[\mathcal F_n=\begin{bmatrix} 
1 & 1 & 1 & \cdots & 1\\
1 & \zeta_{2^n} & \zeta_{2^n}^2 & \cdots & \zeta_{2^n}^{2^n-1}\\
1 & \zeta_{2^n}^2 & \zeta_{2^n}^4 & \cdots & \zeta_{2^n}^{2(2^n-1)}\\
1 & \zeta_{2^n}^3 & \zeta_{2^n}^6 & \cdots & \zeta_{2^n}^{3(2^n-1)}\\
\vdots & \vdots & \vdots &  & \vdots \\
1 & \zeta_{2^n}^{2^n-1} & \zeta_{2^n}^{2(2^n-1)} & \cdots & \zeta_{2^n}^{(2^n-1)^2}\\
\end{bmatrix}.\] For 1-qubit, \[\mathcal F_1= \frac{1}{\sqrt{2}} \begin{bmatrix} 1 & 1 \\ 1 & \zeta_{2} \end{bmatrix}=\frac{1}{\sqrt{2}} \begin{bmatrix} 1 & 1 \\ 1 & e^{\pi i} \end{bmatrix}=\frac{1}{\sqrt{2}} \begin{bmatrix} 1 & 1 \\ 1 & -1 \end{bmatrix}=H.\]i.e. the Hadamard gate is the quantum discrete Fourier  transformation on 1-qubit.

Quantum discrete Fourier transformation $\mathcal F_n$ can be performed by Hadamard gates and CNOT gates.

\section{Deutsch's algorithm}

Deutsch's algorithm is  a simple example of quantum algorithm that shows  computational profit of quantum algorithm. It solves the following problem\footnote{Deutsch, D. Quantum Theory, the Church-Turing Principle and the Universal Quantum Computer. Proceedings of the Royal Society of London A. 400 (1818): 97–117. 1985}.

Let $f:\{0,1\}\longrightarrow \{0,1\}$. If we want to know $f$ is a constant function or not by calculation, we need to calculate $f(0)$ and $f(1)$ classically. However, assuming a quantum gate\[U_f(b_j\otimes b_i):=b_j\otimes b_{f(j)\oplus i}.\] Then, \[(H\otimes I)U_f(H\otimes H)(I\otimes X)(b_0\otimes b_0)=\frac{1}{2}[(1 +(-1)^{f(0)\oplus f(1)})b_0 + (1-(-1)^{f(0)\oplus f(1)})b_1].\] If $f$ is constant, then $f(0)\oplus f(1)=0$ i.e. $P(b_0)=1$. Otherwise, $f(0)\oplus f(1)=1$ i.e. $P(b_1)=1$. In this algorithm, we used $U_f$ only once. 

\section{Shor's algorithm}
RSA is one of the most popular public-key crypto-system. The security of RSA relies on the difficulty of prime factorization. It uses very large two primes $p,q$. The product $pq$ is announced to public and any one who want to sent cryptogram uses $pq$ to encrypt the message. To  Decrypt the cryptogram, one should know what is $p$ and $q$. Since prime factorization is very difficult, one can not find $p$ and $q$ from $pq$.

However, a mathematician Peter Shor published his paper "Algorithms for quantum computation: discrete logarithms and factoring" in 1994 which shows that prime factorization can be obtained fast by his quantum algorithm\footnote{Shor, P.W. Algorithms for quantum computation: discrete logarithms and factoring. Proceedings 35th Annual Symposium on Foundations of Computer Science. IEEE Comput. Soc. Press: 124–134. 1994}.

Let $N$ be the product of two or more odd primes. If we found an element $g\in \mathbb{Z}^*_N$  with the order $|g|$ is even, \footnote{It is proved that there is enough number of element with even order in $ \mathbb{Z}^*_N$.}  $N$ divides $(g^{r/2}-1)(g^{r/2}+1)$ since \[g^r-1\equiv(g^{r/2}-1)(g^{r/2}+1)\equiv0 \mod N.\]  Then $gcd(g^{r/2}-1,N)$ and $gcd(g^{r/2}+1,N)$ are non trivial divisor of $N$. The Euclidean algorithm finding $gcd$ is very fast.

We set two parts of qubit: one part is $u$ which is $n$-qubit system where $N^2\leq 2^n <2N^2$. the other part is $m$-qubit system where $m=\lceil\ln N/\ln 2\rceil $. We operate Shor's algorithm on the $(n+m)$-qubit system $v\otimes u$ as follows.

1. $H$ acts on each qubit in $n$-qubit system $u$:\[(H^{\otimes n}\otimes I^{\otimes m})v_0\otimes u_0=\frac{1}{\sqrt {2^n}}\sum^{2^n-1}_{j=0}v_j\otimes u_0.\]

2. Let $U_x(v_j\otimes u_t):=v_j\otimes u_{t+x^j\mod N}$ for $x\in \mathbb{Z}^*_N$. $U_x$ acts on the $(n+m)$-qubit system:\[U_x\big[\frac{1}{\sqrt {2^n}}\sum^{2^n-1}_{j=0}v_j\otimes u_0\big]=\frac{1}{\sqrt {2^n}}\sum^{2^n-1}_{j=0}v_j\otimes u_{x^j\mod N}.\]

3.  $\mathcal F_n$ acts on the $n$-qubit system: \[ \mathcal F\otimes I^{\otimes m}\big[\frac{1}{\sqrt {2^n}}\sum^{2^n-1}_{j=0}v_j\otimes u_{x^j\mod N}\big]=\frac{1}{2^n}\sum^{2^n-1}_{j=0}\big(\sum^{2^n-1}_{c=0}e^{2\pi ijc/2^n}v_c\big)\otimes u_{x^j\mod N}.\]

4. Carry out the measurement of result \[(\mathcal F\otimes I^{\otimes m})U_x(H^{\otimes n}\otimes I^{\otimes m})(v\otimes u).\] The value of the measurement means the order of $x$.

5. Factorize $N$ using the order of $x$.
\\

 Let the order of $x$ be $r$ and $j=j_0+rk$. $j_0\equiv j \mod r$. $P(v_c\otimes v_{x^{j_0}})$ is \[\frac{1}{2^{2n}}\bigg|e^{2\pi ij_oc/2^n}\sum^{\lfloor2^n/r\rfloor+\delta}_{k=0}e^{2\pi irkc/2^n}\bigg|^2\] where $\delta=0$ ot 1. 

If $r | 2^n$, $P(v_c\otimes v_{x^{j_0}})>0$ only if $2^n/r | c$ and $P(v_c\otimes v_{x^{j_0}})=0$ otherwise. Therefore, the only possible result of measurement is  $v_{t2^n/r}\otimes v_{x^{j_0}}$ for $t \in \mathbb Z$. So, one can find easily the order of $x$.

If $r$ does not divide  $2^n$, one should take a little different process. But the above algorithm still needed.

\section{References}

$\;\;\;\;$[1] D. Coppersmith, An approximate Fourier transform useful in quantum factoring. Technical Report RC19642, IBM. 1994.
\\

[2] D. Deutsch, Quantum Theory, the Church-Turing Principle and the Universal Quantum Computer. Proceedings of the Royal Society of London A. 400 (1818): 97–117. 1985.
\\

[1] G. Mackey, Mathematical Foundations of Quantum Mechanics, W. A. Benjamin, 1963.
\\

[2] J. von Neumann, Mathematical Foundations of Quantum Mechanics, 1932.
\\

[3] P. W. Shor, Introduction to quantum algorithms, arXiv:quant-ph/0005003v2, 2001. 
\\

[4] P.W. Shor,  Algorithms for quantum computation: discrete logarithms and factoring. Proceedings 35th Annual Symposium on Foundations of Computer Science. IEEE Comput. Soc. Press: 124–134. 1994
\\

[5] H. Weyl, The Theory of Groups and Quantum Mechanics, Dover Publications, 1950.

\end{document}